\documentclass[reqno]{amsart}            
\usepackage{amsthm,amssymb,amsfonts,graphicx,cite,color}
\usepackage[bookmarks,bookmarksnumbered,bookmarksopen,colorlinks,backref,linkcolor=blue,citecolor=red]{hyperref}%
\usepackage{xcolor}
\usepackage{graphicx}
\usepackage{mathrsfs}

\usepackage{tikz}
\usetikzlibrary{matrix}

\newtheorem{theorem}{Theorem}[section]
\newtheorem{definition}{Definition}[section]

\newtheorem{corollary}{Corollary}[section]
\newtheorem{remark}{Remark}[section]

\numberwithin{equation}{section}

\begin{document}
\hyphenation{ap-pro-xi-ma-tion}
\title[A mathematical model of the atherosclerosis development in thin blood vessels]
{A mathematical model of the atherosclerosis development in thin blood vessels and its asymptotic approximation}
\author[T.A.~Mel'nyk]{ Taras A. Mel'nyk}
\address{\hskip-12pt  Faculty of Mathematics and Mechanics, Department of Mathematical Physics\\
Taras Shevchenko National University of Kyiv\\
Volodymyrska str. 64,\ 01601 Kyiv,  \ Ukraine
}
\email{melnyk@imath.kiev.ua}

\begin{abstract}
Some existing models of the atherosclerosis development are discussed and
a new improved mathematical model, which takes into account new experimental results about diverse roles
of macrophages in atherosclerosis,  is proposed.
Using technic of upper and lower solutions, the existence and uniqueness of its positive solution are justified.
 After the nondimensionalisation, small parameters are found.
 Then asymptotic approximation for the solution is constructed and justified with the help of asymptotic methods for boundary-value problems
 in thin domains. The results argue for the possibility to replace the complex $3D$ (dimensional) mathematical model with the
corresponding simpler $2D$  model with sufficient accuracy measured by these small parameters.
\end{abstract}

\keywords{Modeling of atherosclerosis; reaction-diffusion system; thin domain; asymptotic approximation
\\
\hspace*{9pt} {\it MOS subject classification:} \  35B40, 92C20, 35K57, 35K50,  74K10
}

\maketitle
\tableofcontents
\section{Introduction}	
Cardiovascular diseases occupy a leading place in mortality in the world. The main cause of these diseases is atherosclerosis.
Therefore the development of atherosclerosis is intensively investigated in last time. There are several theories about atherosclerosis (e.g.~see Refs.~\cite{Guyton,Fan-Wat,Osterud-Bjorklid,Ross}), but none can explain the whole process because this disease is associated with many risk factors. The atherosclerotic process is not fully understood till now.
However, many researchers agree that the damage or dysfunction of the arterial endothelium  and high level of low
density lipoproteins in blood vessels play the main role in the development of atherosclerosis.

In existing mathematical models, researchers are trying to take into account various factors and different types of molecules involved in the development of this illness. As a result we have models containing from two to more then thirty differential equations.
However, even a such cumbersome system cannot take into account all features of atherosclerosis.

Starting from a review of current modeling approaches of atherosclerosis, this article aims at specifying the optimal prospects for research on the mathematical study of atherosclerosis involving special rigorous asymptotic method that enables the reasonable approximation of the original model with its features. The main idea is that atherosclerosis should be looked at as complex system of enzyme reactions that are able to choose their individual dynamics.

The paper is organized as follows.
Section~\ref{Sec2} describes the basic stages of  mechanism of the atherosclerosis development. Moreover, it presents some preliminary
models focusing on their shortcomings. Taking into account those shortcomings, Section~\ref{Sec3}  offers a new mathematical model of atherosclerosis. Using technic of upper-lower solution,  the existence and uniqueness of its positive solution are justified. After the nondimensionalisation, we find small parameters  of the mathematical model and make the asymptotic analysis as those small parameters tend to zero in Section~~\ref{Sec4}. Namely, we find the corresponding limit problem, construct the asymptotic approximation, find its residuals, estimate them and prove asymptotic estimates for the difference between the solution and the approximating function. In Section~\ref{Sec5}, we discuss several generalizations and research perspectives.

\section{Existing mathematical models of the atherosclerosis development}\label{Sec2}

For convenience of readers we present the following short glossary:
\begin{itemize}
\item
{\bf LDLs} and {\bf HDLs}  --  low and high density lipoproteins; they transport lipid and cholesterol
to cells and tissues around the body;
\item
{\bf Free radicals} --  extremely reactive molecules which scavenge electrons causing oxidation of the target;
\item
{\bf Monocytes} --  the largest of the white blood cells (leukocytes); they are able to get in  places of
inflammation or tissue damage;
\item
{\bf Cytokines} -- small proteins that are important in cell signaling; they are secreted by certain cells
of the immune system and have an effect on other cells including the stimulation or inhibition of their growth
 and other functional activities;
\item
{\bf Intima} -- the first thin layer of the blood vessel wall after the endothelium.
\end{itemize}

The mechanism of the atherosclerosis development can be shortly sketched as follows
(see e.g. Ref.~\cite{Fan-Wat}):

\centerline{{\bf Stage 1}}
\begin{enumerate}
\item[{\bf 1.}]
After damage of the vessel endothelium, circulating lipids (mostly (LDLs)) begin to accumulate
at the site of the injury.
  \item[{\bf 2.}]
  The process of atherosclerosis begins when LDLs penetrate into the intima of
the arterial wall where they are oxidized (ox-LDLs).
\end{enumerate}

\centerline{{\bf Stage 2}}

\begin{enumerate}
  \item[{\bf 3.}]
  Ox-LDLs in the arterial intima is considered by the immune system
as a dangerous substance, hence an immune response is launched:
monocytes circulating in the blood adhere to the endothelium and then they penetrate to the intima.
  \item[{\bf 4.}]
  Once in the intima, these monocytes are converted into macrophages.
\end{enumerate}

\centerline{{\bf Stage 3}}

\begin{enumerate}
  \item[{\bf 5.}]
The macrophages phagocytose (ingest) the ox-LDLs and become foam cells. Simultaneously
chronic inflammatory reaction is started:
\begin{itemize}
  \item
  macrophages secrete pro-inflammatory cytokines that promote the recruitment of new monocytes and thereby
  support the production of new pro-inflammatory cytokines;
  \item
  macrophages with a large intake of the ox-LDLs become foam
cells that cause the growth of the intimal layer and thereby amplify the endothelial dysfunction;
\item
 this auto-amplification phenomenon is compensated by an anti-inflammatory response mediated by anti-inflammatory
cytokines. 
\end{itemize}
\end{enumerate}
This is the short sketch of the atherosclerosis development.
In reality, many complex biochemical reactions are hidden behind of these steps. Some of them will be discussed later.

Now let us consider existing mathematical models.
The first models of penetration of cholesterol in the arterial wall (the first step in the stage~1) were considered in Refs.~\cite{BraColSmith,many1990}.
These models are one dimensional systems of differential equations with respect to the radial coordinate
(the vessel wall was assumed to be uniform along its length).

In Ref.~\cite{KhaGenKazVol-07} the authors have developed  simple models of the reactions arising in the arterial intima.
 The first model represents the following one-dimensional reaction-diffusion system:
\begin{equation}\label{sys1}
\left\{
    \begin{array}{rcl}
\partial_{t}u_m  &=& d_m \partial^2_{xx} u_m + f_1(u_c) - \beta u_m,
                                        \\[2mm]
\partial_{t}u_c  &=& d_c \partial^2_{xx} u_c + f_2(u_c) u_m  - \gamma u_c + b,
\end{array}
\right.
\end{equation}
for $x \in (0, L).$  Here $\partial_{t}u = \frac{\partial u}{\partial t},$  $\partial^2_{xx}u = \frac{\partial^2 u}{\partial x^2},$
the value $u_m$ is  the concentration of monocytes, macrophages and foam cells together  in the intima, $u_c$ is the concentration of cytokines.
The function $f_1(u_c)$  describes the recruitment of monocytes from the blood flow, $f_2(u_c) u_m$ is the rate of production of the cytokines
which depends on their concentration and on the concentration of the blood cells:
\begin{equation}\label{fun1}
f_1(u_c) = \frac{\alpha_1 + \beta_1\, u_c}{1 + u_c/ \tau_1}, \qquad  f_2(u_m) = \frac{ \alpha_2\, u_m}{1 + u_m / \tau_2}.
\end{equation}
The negative terms ``$- \beta u_m$'' and
``$- \gamma u_c$'' correspond to the natural death of the cells and chemical substances, and
the last term $b$ in the right-hand side describes the ground level of the cytokines in the intima.

The second model  deals  with  the system
\begin{equation}\label{sys2}
\left\{
    \begin{array}{rcl}
\partial_{t}u_m  &=& d_m \Delta_{x_1x_2} u_m  - \beta u_m,
                                        \\[2mm]
\partial_{t}u_c  &=& d_c \Delta_{x_1x_2} u_c + f_2(u_c) u_m  - \gamma u_c + b,
\end{array}
\right.
\end{equation}
in the two-dimensional thin rectangle
$$
\Omega_\varepsilon = \{(x_1, x_2)\in \Bbb R^2: \ \ x_1 \in (0, L), \ \  x_2 \in (0, \varepsilon)\}
$$
(here $\varepsilon$ is a small parameter, \ $\Delta_{x_1x_2} u := \partial^2_{x_1x_1} u + \partial^2_{x_2x_2} u$) with
the nonlinear boundary condition that takes  into account the recruitment of monocytes from the blood flow:
\begin{equation}\label{sys3}
   d_m \,\partial_{x_2} u_m = \varepsilon\, f_1(u_c)  \quad \text{at} \ \ x_2=\varepsilon, \ \ x_1 \in (0, L),
\end{equation}
and the homogeneous Neumann boundary condition on the rest part of the boundary $\partial \Omega_\varepsilon;$
the homogeneous Neumann boundary conditions are imposed everywhere
on $\partial \Omega_\varepsilon$ for  the cytokine concentration $u_c$.

The authors have analyzed one-dimensional model depending on the parameters $\alpha_1, \alpha_2, \beta_1, \tau_1, \tau_2$, proved the existence of travelling wave solutions and with the help formal asymptotic analysis  (see the appendix~B in Ref.~\cite{KhaGenKazVol-07}) showed that the system (\ref{sys1}) can obtained from (\ref{sys2})--(\ref{sys3}) as $\varepsilon \to 0.$

In Ref.~\cite{KhaGenKazVol-12} the authors continued to study two-dimensional system (\ref{sys2}) in the strip
 $\Omega_h = \{(x_1, x_2)\in \Bbb R^2: \ \ -\infty <x_1 < +\infty, \ \ 0< x_2 < h\}$ with
the following boundary conditions:
\begin{gather}
\partial_{x_2} u_m = 0, \ \ \partial_{x_2} u_c = 0  \quad \text{at} \ \ x_2=0, \notag
\\
\partial_{x_2} u_m = f_1(u_c), \ \ \partial_{x_2} u_c = 0 \quad \text{at} \ \ x_2=h, \label{bc-1}
\end{gather}
where the functions $f_1$ and $f_2$ are sufficiently smooth and satisfy the following conditions:
\begin{gather}\label{cond0}
f_2(u) > 0 \ \ \text{for} \ \ u > 0, \quad f_2(0)=0, \quad f_2(u) \to f_2^+ \ \ \text{as} \ \ u \to +\infty;
\\
f_1(u) > 0 \ \ \text{for} \ \ u > c_0,  \quad f_1(c_0)=0, \quad f_1(u) \to f_1^+ \ \ \text{as} \ \ u \to +\infty, \label{cond0-1}
\end{gather}
and $f_1'(u) > 0.$ The positivity of the solution of the system (\ref{sys2})--(\ref{sys3}) and
the existence of travelling waves in the monostable case are proved. Results of numerical simulations for those systems were also obtained
in Refs.~\cite{KhaGenKazVol-07,KhaGenKazVol-12}.

These models include only two variables $u_m$ and $u_c$ (it is essentially for the method proposed by the authors). Obviously, it is not enough to describe the full picture of atherosclerosis. In addition, these models omit some essential features of atherosclerosis development, namely,
 LDLs penetration and oxidation (see   the stage 1); \
  transformation of monocytes into macrophages and then into foam cells,
 the chronic inflammatory reaction (see   the stage 3);
\ the diffusion coefficients of macrophages and cytokines are quite different
$(8.64\times10^{-7}$ and $1.08\times 10^2$ $cm^2\, day^{-1},$ see e.g. Ref.~\cite{HaoFried-14})
 and  cytokines are active in very low concentrations and their secretion is short and strictly regulated
(see e.g. Refs.~\cite{AitTalMalTed,cytokines-book}). Also,  the size of the damaged endothelium
and cylindrical shape of vessels were not taken into account (boundary conditions (\ref{sys3}) and  (\ref{bc-1}) are set on the whole side both of the rectangle and strip).

At first time, all three stages and in addition the formation of a plaque were mathematically simulated in Ref.~\cite{HaoFried-14}.
The model includes the following key variables:
LDLs and HDLs, free radicals and ox-LDLs, \ six types of cytokines (MCP-1, IFN-$\gamma$, IL-12, PDGF,  MMP and TIMP),
and four types of cells (macrophages, foam cells, T-cells and smooth muscle cells). This model consists of eighteen partial differential equations in a  plane domain. Unfortunately,  the existence and uniqueness of the solution of this system were not justified,  the nondimensionalization was not made, monocytes were not involved in the model. In addition,  there are many special assumptions, e.g.:  \ {\it all cells are of the same volume and surface area, so that the diffusion coefficients of the all cells have the same coefficient; the cells are moving with a common velocity, \ the domain has polygonal boundaries}.

But, even those equations are not enough to describe a complete model of atherosclerosis.
For instance, in Refs.~\cite{BujaMcAllister-07,UsmRibSadGil-15} the authors showed that  more than twenty substances (different cells, various cytokines, hormones, chemical mediators and effectors) are involved in the pathogenesis of atherosclerosis.

Now researchers study different stages of atherosclerosis in more detail.
In Ref.~\cite{CobSheMax}  the preliminary lipoprorein oxidation model includes nine ordinary equations, the extended model has twelve equations.
The paper\cite{GuiShiSunAkaMur-12} discusses the central roles of macrophages in different states of atherosclerosis, focusing on the role of inflammatory biomarkers in predicting primary cardiovascular events. There (Section 4) it was showed how many different
pro- and anti-inflammatory cytokines, chemokines, mediators, enzymes and biomarkers related only to macrophages
are involved in the progression of atherosclerosis.

So, the following  question arises. How  should we mathematically simulate the atherosclerosis development?
Should we collect all those equations together? Then we get a system with about forty equations.

\section{Improved mathematical model}\label{Sec3}

It is known that many disease are a set of biochemical reactions.
Most of these reactions are enzymatic reactions that take into account the influence of smaller molecules  (like cytokines and enzymes).
Therefore, the answer on the question above is the balance between biological  realism and simplicity which can be described by
known scenarios of enzymatic reactions and their mathematical models.

The mechanism of the most basic enzymatic reaction (a substance S reacting with an enzyme E to produce a substance P),
first proposed by Michaelis and Menten\cite{MichMenten}, is represented schematically by
$$
S + E \ \ {\rightleftarrows}^{k_1}_{k_{-1}} \ SE \ \ \stackrel{k_2}{\longrightarrow}\ \ P + E
$$
and described by the system of differential equations (see e.g. Ref.~\cite{Murray}\,[Sec.~6.1])
\begin{gather*}
\frac{ds}{dt}= - k_1  e \, s + k_{-1} c , \quad \frac{de}{dt}= - k_1  e\, s + (k_{-1} + k_2) c ,
\\
\frac{dc}{dt}=  k_1  e\, s - (k_{-1} + k_2) c , \quad \frac{dp}{dt}=  k_2 c,
\end{gather*}
with the initial conditions: $s(0)=s_0, \ e(0)=e_0, \ c(0)=0, \ p(0)=0.$
Here $s=[S],$ $e=[E],$ $c=[SE],$ $p=[P]$ are corresponding concentrations,  $k$'s are the rate constants.

Taking into account the second and third initial conditions, we can reduce this system to the following two differential equations:
$$
\frac{ds}{dt} = - f(s),\qquad \frac{dp}{dt} = f(s),
$$
where $\displaystyle{f(s)=\frac{k_2 \, e_0 \, s}{K_m + s}}$ and   $\displaystyle{K_m=\frac{k_{-1}+ k_2}{k_1}}$ is the Michaelis constant.

Thus, we don't need any differential equation for the enzyme $E$
to describe the enzymatic reaction that  converts the substance $S$ to $P.$

Typical choices of such functions in biological applications are
\begin{equation}\label{enzym-functions-1}
  f(s)= \frac{\lambda \, s^p}{\mu + s^p} \quad \text{or} \quad f(s)= \frac{\alpha + \lambda \, s^p}{\mu + s^p} \quad (p > 0,  \ \ \mu > \frac{\alpha}{\lambda}, \ \ \lambda >0);
\end{equation}
and more general case
\begin{equation}\label{enzym-functions-2}
f(s)= \frac{\lambda \, g(s)}{\mu + g(s)}   \quad (g'(s) \ge 0, \ \ g(0)=0, \ \ \mu >0).
\end{equation}

In the present paper this approach will be used to describe biochemical reactions in the development of atherosclerosis.

From the review above, we see that the following molecules and cells
\begin{itemize}
\item
LDLs $(L),$  ox-LDLs $(L_{ox}),$
monocytes $(m),$  macrophages $(M),$  foam cells $(F)$
\end{itemize}
qualitatively reproduce the atherosclerosis development.
Therefore, we consider their concentrations $L, L_{ox}, m, M, F$ as key variables in our model and do not include other molecules with very low concentration.

\begin{figure}[htbp]
\centering
\includegraphics[width=6cm]{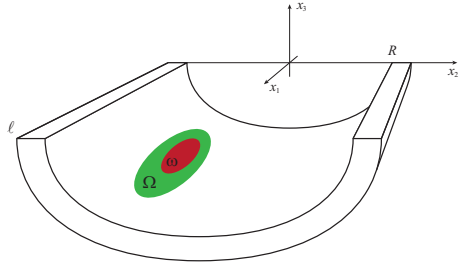}
\caption{The longitudinal cross-section of the domain $C_{R,\rho_0}$}\label{f1}
\end{figure}

Now let us go to modelling. At first, we define the following pipe domain (see Fig.~\ref{f1}):
$$
C_{R,\rho_0} := \Big\{x=(x_1, x_2, x_3) \in \Bbb R^3: \ \ x_1 \in (0, l ) , \quad \sqrt{x_2^2+x_3^2}\, =: \rho \in (R, R+\rho_0) \Big\}
$$
that will be a prototype of an artery. Its boundary consists of two bases
$$
  \Gamma_{0,l}:= \partial C_{R,\rho} \cap \big(\{x: \  x_1=0\}\cup \{x: \  x_1= l\}\big)
$$
and two cylindrical parts
$$
  \Upsilon:= \partial C_{R,\rho} \setminus  \Gamma_{0,l}.
$$
In the inner cylindrical part $\Upsilon_R :=  \big\{x\in \Bbb R^3: \ \ x_1 \in (0, l) , \  \rho = R \big\},$ we consider two
smooth surfaces $\omega$ and $\Omega$ such that
$
\omega \subset \Omega \ \subseteq \ \Upsilon_R.
$
The surface $\omega$ is a prototype of the damage or dysfunction of the arterial endothelium, and $\Omega$ is a surface of
the penetration of monocytes into the intima.

\subsection{Modelling the penetration of LDLs and monocytes}\label{st-1}

The penetration of LDLs was simulated by usual linear Robin boundary condition (see equastion (23) in Ref.~\cite{HaoFried-14}).
However, as has been found in Ref. \cite{Fan-Wat}, the damaged endothelium is involved in the activation of endothelial adhesion molecules: \ (ICAM-1),\ (VCAM-1), \ P-selestin, \ and E-selestin that attract LDLs to endothelial cells. Thus, the penetration of LDLs is not pure diffusion.
 We admit that the penetration is obeyed to  the Michaelis-Menten rule (when enzyme reactions appear at the boundary surface of the diffusion medium it leads to nonlinear boundary conditions~\cite{Ross,Pao}. So, we propose the following boundary condition:
\begin{equation}\label{penetration-con}
d_L \partial_{\nu} L =  \alpha_0 \frac{L}{k_0 + L} + {\large {\copyright_1}} \quad \text{on} \ \ \omega,
\end{equation}
and $\partial_{\nu} L = 0 $ on $\Upsilon\setminus \omega$ (the non-flux boundary condition); \ $\partial_{\nu} L$ is the outward normal derivative.
 Hereinafter all coefficients in differential equations and boundary conditions are positive (biological meaning of the coefficients $\alpha_0$ and $k_0$ before  (\ref{gamma_1})).   The second term ``${\copyright_1}$'' in (\ref{penetration-con}) simulates  another process and will be presented later.

\medskip

The key to the early inflammatory response is the activation of the endothelial cells through special enzyme (apoB-LPs), which is intensified by ox-LDLs, that leads to the recruitment of monocytes from the lumen into the intima~\cite{MesLey}. Therefore, we model this enzyme reaction by the following boundary condition:
\begin{equation} \label{bc-2-mono}
d_m \, \partial_{\nu} m = \frac{\alpha_1 \, L_{ox}}{k_1 + L_{ox}} + {\large {\copyright_2}} \quad \text{on} \ \ \Omega \supset \omega;
\end{equation}
in addition the non-flux boundary condition on $\Upsilon\setminus \Omega.$ The term ``${\copyright_2}$'' will be presented later.

\subsection{Modelling the efflux of macrophages}\label{st-1+}
Recently, there has been a great deal of interest in the role of macrophages in atherosclerotic disease (see e.g. Refs.~\cite{JohnNewby,ManGarLoc,Paulson-others,ManDiv-11,GuiShiSunAkaMur-12,Feig-others}). It was discovered that monocytes are transformed into two different kind of macrophages $M_1$ and $M_2:$
\begin{equation}\label{diagram}
\begin{array}{rcl}
& & M_1 \ \ (\text{classically activated;} \ \ 85\% \sim 95\%)
\\
 & \nearrow &
\\
m &  &
\\
& \searrow &
\\
& & M_2 \ \ (\text{alternatively activated})
\end{array}
\end{equation}
Classically activated macrophages $(M_1)$ promote the inflammation
and product cytokines $IL$-23, $IL$-6, $IL$-12, $IL$-1, $TNF$-$\alpha;$
while the second type macrophages $(M_2)$ secrete anti-inflammatory cytokines $IL$-4, $IL$-13, $IL$-1, $IL$-10 (see Ref.~\cite{GuiShiSunAkaMur-12}).
In addition, it is turned out that $M_2$-macrophages with the ingested ox-LDLs inside can return (efflux) to the blood flow in the vessel~\cite{Feig-others}. We model this by
\begin{equation}\label{bc-3}
d_{M_2} \partial_{\nu} M_2 = - \frac{\alpha_2 \, M_2^{p_2}}{k^{p_2}_2 + M_2^{p_2}}  \quad \text{on} \ \ \Omega \supset \omega \quad (p_2 > 1);
\end{equation}
and the non-flux boundary condition on $\Upsilon\setminus \Omega.$

\begin{remark}
 Hereinafter,  parameters $p_1, p_2, p_3$ are greater than $1.$  This means some delay in time and less mobility
of macrophages and foam cells because of their size.
\end{remark}

\subsection{Processes in the intima}

Evolution and oxidation of LDLs are modelled by the following two reaction-diffusion equations:
\begin{gather}
\partial_t L - d_L \Delta_x L= - \Theta(x, t) \, L + \psi(x, t),\label{eq-1}
\\[2mm]
\partial_t L_{ox} - d_{L_{ox}} \Delta_x L_{ox}=  \Theta(x, t) L   -  \underbrace{\lambda_{L_{ox}M_1} \frac{L_{ox}}{K_{M_1} + L_{ox}}\, M_1}_{intake \, of \, L_{ox}\, by \, M_1}
-  \underbrace{\lambda_{L_{ox}M_2} \frac{L_{ox}}{K_{M_2} + L_{ox}}\, M_2}_{intake \, of \, L_{ox}\, by \, M_2} \label{eq-2}
\end{gather}
in the domain $C_{R,\rho},$ where $d_L$ and $d_{L_{ox}}$ are the diffusion coefficients of
LDLs and ox-LDLs in the intima respectively. The function $\psi$ is smooth, nonnegative and  its compact support is contained in $\omega$
(its biological meaning in Remark~\ref{rmk-1}).

In papers Refs.~\cite{HaoFried-14,CobSheMax}, the term ``$k_L r L$'' is used to simulate  the production of ox-LDLs
in the reaction of LDLs with free radicals $r$ whose baseline growth is a bifurcation parameter $r_0.$ However, as was noted in Refs.~\cite{CobSheMax}, it is impossible to measure $r_0$ since many social factors such as drinking, smoking and ecological factor
highly and unpredictably increase the parameter~$r_0.$ Therefore, we use a smooth function $\Theta(x, t),$ $(x, t)\in C_{R,\rho_0}\times [0, T],$ that can depend both on the concentration of free radicals $r$  and on other values (like HDLs, vitamins $C$ and $E$) and assume that
$\Theta|_{\Gamma_{0,l}}=0,$ $\Theta|_{t=0}=0,$ and $0\le \Theta \le c_2.$

The last two terms in the right-hand side of equation (\ref{eq-2}) are  reduction terms that represents the phagocitosis (intake) of  ox-LDLs by macrophages~$M_1$ and $M_2$ respectively. By virtue of the diagram (\ref{diagram}) we regard that the reduction rate  $\lambda_{L_{ox}M_2}$
 is much lower than $\lambda_{L_{ox}M_1}$ and  $K_{M_2} < K_{M_1}.$

In Ref.~\cite{HaoFried-14} \  the authors did not take into account two different types of macrophages
and have chosen  the term ``$\lambda_{L_{ox} M} L_{ox} M$'' to simulate the ingestion of  ox-LDLs.
However, we prefer terms like   $\lambda \frac{L_{ox}}{K_M + L_{ox}}\, M$ since they appear  in many modern prey-predator models and other biochemical reactions~\cite{Perthame}.

The equation (\ref{eq-2}) is supplemented with  the boundary condition $\partial_{\nu} L_{ox} = 0$ on~$\Upsilon$.

\medskip

Equations, which model the evolution of ``$m$'',  ``$M_1$'', and ``$M_2$''  in the intima, are as follows:
\begin{equation}\label{eq-3}
\partial_t m - d_m \Delta_x m=  - \underbrace{\lambda_{m M_1} \frac{m}{k_{M_1} + m}}_{transformation \, into \, M_1}
- \underbrace{\lambda_{mM_2} \frac{m}{k_{M_2} + m}}_{transformation \, into \, M_2}
\ \ \underbrace{- \ \ \beta_m\, m}_{death \,of\, m},
\end{equation}
\begin{equation}\label{eq-4}
\partial_t M_1 - d_{M_1} \Delta_x M_1 =  \lambda_{m M_1} \frac{m}{k_{M_1} + m} + \lambda_{L_{ox}M_1} \frac{L_{ox}}{K_{M_1} + L_{ox}}\, M_1
 - \underbrace{\lambda_{M_1 F} \frac{M_1}{K_F + M_1}}_{transformation \, into \, F}\ - \ \beta_{M_1} M_1,
\end{equation}
\begin{equation}\label{eq-5}
\partial_t M_2 - d_{M_2} \Delta_x M_2 =  \lambda_{mM_2} \frac{m}{k_{M_2} + m} + \lambda_{L_{ox}M_2} \frac{L_{ox}}{K_{M_2} + L_{ox}}\, M_2 \
 - \ \beta_{M_2} M_2,
\end{equation}
where  $\lambda_{m M_2} < \lambda_{m M_1}$ and $k_{M_2} < k_{M_1}$ because of (\ref{diagram}).
We would like to note once again that monocytes were not involved in the model of  the paper Ref.~\cite{HaoFried-14}.

The remained macrophages ``$M_1$'' with a large intake of ``$L_{ox}$'' become  foam cells~``$F$''.
This reaction can be simulated by the following equation:
\begin{equation}\label{eq-6}
\partial_t F - d_{F} \Delta_x F =  \lambda_{M_1 F} \frac{M_1}{K_F + M_1} \ - \ \beta_F F
\end{equation}
in the intima. Obviously, $\partial_{\nu} F = \partial_{\nu} M_1 = 0$ on $\Upsilon.$

\begin{remark}
Since the transformation of a substance $S$ into a substance $P$ is a set of enzyme reactions, we use terms like $\lambda \frac{S}{k_P + S}$
in equations (\ref{eq-3}) -- (\ref{eq-6}). However, we can take more general terms (\ref{enzym-functions-1}) and (\ref{enzym-functions-2}) without any restrictions. In  reality, when selecting the view of each term, experimental study should be taken into account.
\end{remark}

\subsection{Modelling of the stage 3 (chronic inflammatory reaction)}

As was mentioned above $M_1$-macrophages secrete pro-inflammatory cytokines
that promote the recruitment of new monocytes; on the other hand, $M_2$-macrophages secrete anti-inflammatory cytokines that inhibit the recruitment of new monocytes. We model this  by the term ``${\copyright_2}$''  in the boundary condition (\ref{bc-2-mono}), namely
\begin{equation}\label{bc-5}
d_m \partial_{\nu} m = \frac{\alpha_1 L_{ox}}{k_1 + L_{ox}} + \frac{\alpha_3 M_1^{p_1}}{(k_3 + M_2)(k_4^{p_1} + M_1^{p_1})} \quad \text{on} \ \ \Omega \quad (p_1 >1).
\end{equation}

\begin{remark}
The activator-inhibitor system was firstly considered by Gierer and Meinhardt~\cite{GerMei}
in a pattern formation model (the $``u``$-activator -- $``v``$-inhibitor term in this model is
${\frac{u^2}{v(1+ K u^2)}})$
and now it has been used in many applications (see e.g. Ref.~\cite{Murray-2}).
\end{remark}

The foam cells cause the growth of the intimal layer, thereby they amplify the endothelial dysfunction
and additional penetration of LDLs. This feedback can be simulated by the term ``${\copyright_1}$''  in the boundary condition (\ref{penetration-con}), namely
\begin{equation}\label{bc-4}
d_L \partial_{\nu} L =  \alpha_0\, \frac{L}{k_0 + L} + \frac{\alpha_5 F^{p_3}}{k_5^{p_3} + F^{p_3}}  \quad \text{on} \ \ \omega \quad (p_3  >1).
\end{equation}

\begin{remark}
Such kind of feedback loop is common in the biochemical control circuit of substances $u_1,\ldots, u_n,$ which  is described by the system of equations
\begin{gather*}
\frac{d u_1}{dt} = g(u_n) - k_1 u_1,
\\
\frac{d u_i}{dt} = u_{i-1} - k_i u_i, \quad i=2,\ldots,n \quad (k_i >0),
\end{gather*}
where $g$ is the feedback function. Typical choice of $g$ in applications is
$$
g(u)= \frac{\nu + \alpha u^p}{K + u^p} \quad (K > 0, \ \ \nu, \alpha \ge 0, \ \ p > 0).
$$
The first model of this type was proposed by Goodwin~\cite{Goodwin} with the following feedback function:
$
\displaystyle{g_G(u) = \frac{\nu_0}{1+ (\frac{u}{k})^p}}.
$
\end{remark}

\subsection{Nondimensionalization}

The equations (\ref{eq-1}), (\ref{eq-2}), (\ref{eq-3})--(\ref{eq-6})  with the boundary conditions mentioned above,
with the Dirichlet  conditions
$$
L = L_{ox} = m = M_1 = M_2 = F = 0 \ \ \text{on} \ \Gamma_{0, l},
$$
and with the initial conditions $L= L_{ox}=  m= M_1=0 = M_2= F=0$ at $t=0$
form the improved mathematical model of atherosclerosis development.

\begin{remark}\label{rmk-1}
It is easy to verify that the compatibility condition at  $t=0$ is satisfied. General case  $L|_{t=0}=L_0,$ which describes the LDL concentration at the beginning of the atherosclerosis development, can be reduced with the substitution $L - L_0$ to the homogeneous one. This substitution explains
the appearance of the function $\psi$ in (\ref{eq-1}) that  depends on the LDL concentration in the blood.
\end{remark}

Now we nondimensionalise the equations by setting
$$
u_1= \frac{L}{k_0}, \ \ u_2= \frac{L_{ox}}{k_1},  \ \ u_3= \frac{m}{k_{M_1}}, \ \ u_4 = \frac{M_2}{k_2}, \ \ u_5= \frac{M_1}{k_4},
\ \ u_6= \frac{F}{k_5},
$$
$$
t^*= \frac{d_L}{R^2}\, t, \quad x^*= \frac{x}{R}.
$$

Then the pipe domain $C_{R, \rho_0}$ is transformed into a new pipe domain
$$
C_{\varepsilon} := \Big\{x^* \in \Bbb R^3: \ \ x^*_1 \in (0, \ell ) , \quad r^* \in (1, 1+\varepsilon) \Big\}, \quad \ell=\frac{l}{R}, \ \
\varepsilon=\frac{\rho}{R}.
$$
Its boundary consists of the bases (we drop the asterisks for algebraic simplicity)
$$
  \Gamma^\varepsilon_{0,\ell}:= \partial C_{\varepsilon} \cap \big(\{x: \ \ x_1=0\}\cup \{x: \ \ x_1=\ell\}\big)
$$
and two cylindrical parts
$
  \Upsilon_{1, \varepsilon} := \Upsilon_1 \cup \Upsilon_{1+\varepsilon}.
$
The surfaces $\omega$ and $\Omega$ are transformed into two cylindrical smooth surfaces $\omega_1$ and $\Omega_1$ respectively that belong to $\Upsilon_1.$

Let us denote by ${\bf u} := (u_1, \ldots, u_6).$ Then the nondimensional system is as follows:
\begin{equation}\label{start_prob}
\left\{
    \begin{array}{rcll}
    \partial_t {\bf u} - \mathfrak{D} \Delta_x {\bf u}  &=& {\bf F}({\bf u})
                    &\quad \mbox{in} \ C_{\varepsilon}\times (0,T),
                    \\[2mm]
 \mathfrak{D}  \partial_{\nu} {\bf u} &=&   {\bf G}({\bf u})
                    &\quad \mbox{on} \ \Upsilon_{1}\times (0,T),
                    \\[2mm]
 \mathfrak{D}  \partial_{\nu} {\bf u} &=&   {\bf 0}
                    &\quad \mbox{on} \ \Upsilon_{1+\varepsilon}\times (0,T)
                    \\[2mm]
  {\bf u} &=& {\bf 0}&\quad \text{on} \  \Gamma^\varepsilon_{0, \ell}\times (0,T),
 \\[2mm]
    \left. {\bf u}\right|_{t=0} &=&{\bf 0} &\quad \mbox{in} \   C_{\varepsilon},
   \end{array}\right.
\end{equation}
where  the matrix $\mathfrak{D} = {\rm diag}\left(1,d_2,\ldots,d_6\right)$
introduces  the diffusion  constants
$$
d_2=\frac{d_{L_{ox}}}{d_L}, \quad d_3=\frac{d_{m}}{d_L}, \quad d_4=\frac{d_{M_2}}{d_L},\quad d_5=\frac{d_{M_1}}{d_L},\quad d_6=\frac{d_{F}}{d_L};
$$
the reaction terms
\begin{equation}\label{reaction-term-1}
{\bf F}({\bf u}) = 
\left(
\begin{array}{l}
f_1({\bf u}) = - \, \Theta(x,t) \, u_1 +  \psi(x,t)
\\[2mm]
f_2({\bf u}) = \mu_1\, \Theta(x,t) \, u_1 - \mu_2 \dfrac{ u_2}{\lambda_1 + u_2} \, u_5 -  \mu_3 \dfrac{ u_2}{\lambda_2 + u_2} \, u_4
\\[3mm]
f_3({\bf u}) = -  \mu_4 \dfrac{ u_3}{1 + u_3}   - \mu_5 \dfrac{ u_3}{\lambda_3 + u_3} -  \delta_3 \, u_3
\\[3mm]
f_4({\bf u}) = \mu_6 \dfrac{ u_3}{1 + u_3}  + \mu_7 \dfrac{u_2}{\lambda_4 + u_2} \, u_4  -  \delta_4 \, u_4
\\[3mm]
f_5({\bf u}) = \mu_8 \dfrac{ u_3}{\lambda_5 + u_3}  +  \mu_9 \dfrac{u_2}{\lambda_6 + u_2} \, u_5 - \mu_{10} \dfrac{ u_5}{\lambda_7 + u_5}  - \delta_5 \, u_5
\\[3mm]
f_6({\bf u}) = \mu_{11} \dfrac{ u_5}{\lambda_8 + u_5}  - \delta_6 \, u_6
\end{array}
\right)
\end{equation}
and
\begin{equation}\label{reaction-term}
{\bf G}({\bf u}) =
\left(
\begin{array}{l}
g_1({\bf u}) = \varphi_1(x)\left(
\gamma_1 \dfrac{u_1}{1 + u_1} + \gamma_2 \dfrac{ u_6^{p_3}}{1 + u_6^{p_3}}\right)
\\
g_2({\bf u}) = 0
\\
g_3({\bf u}) = \varphi_2(x)\left(
\gamma_3 \dfrac{ u_2}{1 + u_2} + \gamma_4 \dfrac{ u_5^{p_1}}{(\lambda_{9} + u_4)(1 + u_5^{p_1})}\right)
\\[4mm]
g_4({\bf u}) = - \varphi_2(x) \gamma_5 \dfrac{  u_4^{p_2}}{1 + u_4^{p_2}}
\\
g_5({\bf u}) = 0
\\
g_6({\bf u}) = 0
\end{array}
\right).
\end{equation}
Since the boundary conditions (\ref{bc-3}), (\ref{bc-4}) and (\ref{bc-5}) are localized on $\omega$ and $\Omega$ respectively,
in (\ref{reaction-term}) we introduce special smooth functions $\varphi_1$ and $\varphi_2$ such that
$\text{supp}(\varphi_1) \subset \omega_1,$ $ 0< \varphi_1 \le 1$ on $\omega_1,$
$\text{supp}(\varphi_2) \subset \Omega_1,$ $0< \varphi_2 \le 1$ on $\Omega_1.$
All constants $\{\lambda_i\},$ $\{\mu_i\},$ $\{\delta_i\},$ and $\{\gamma_i\}$ are positive; $\mu_3 < \mu_2,$ $\lambda_2 < \lambda_1,$
$\mu_5 < \mu_4,$ $\mu_8 < \mu_6,$ $\lambda_3 < 1,$ $\lambda_5 < 1.$

\subsection{Existence and uniqueness of the solution to problem (\ref{start_prob})}\label{subs-exist}

For this we will use the method of upper and lower solutions, which was developed in Ref.~\cite{Pao}. This method leads not only to the basic results of existence and uniqueness of solutions but also to some their qualitative properties. For the convenience of readers we present some definitions from this book  in adapting to problem (\ref{start_prob}).

\begin{definition}
A vector-function ${\bf f}=(f_1,\ldots,f_6)$ $(f_i\equiv f_i(t,x,u_1,\ldots,u_6))$ is said to possess a quasi-monotone property if for any
$ i\in\{1,\ldots,6\}$ there exist $ a_i, b_i \in \Bbb N_0, \ a_i + b_i = 5,$ such that
 $f_i(\cdot,\cdot,u_i, [{\bf u}_{a_i}], [{\bf u}_{b_i}])$ is monotone nondecreasing in   $[{\bf u}_{a_i}]$ and
is monotone nonincreasing in $[{\bf u}_{b_i}].$ Here $[{\bf u}_{a_i}]$, $[{\bf u}_{b_i}]$ denote the $a_i$-components and $b_i$-components
of the vector ${\bf u} \in \Bbb R^6_+,$ respectively.
\end{definition}

For example,  the function
$$
f_2= \mu_1 \Theta(x,t)\, u_1 - \mu_2 \frac{ u_2}{\lambda_1 + u_2} \, u_5 -  \mu_3 \frac{u_2}{\lambda_2 + u_2} \, u_4
$$
is monotone nondecreasing in $u_1,$ $u_3,$ $u_6$ and is
monotone nonincreasing in $u_4,$ $u_5$ (recall that the function $\Theta$ is nonnegative (see subsection~\ref{st-1})); thus
$a_2=3,$ $b_2=2.$

It is easy to verify that the vector-functions ${\bf F}$ and ${\bf G}$ are quasi-monotone.

\begin{definition}
A pair of functions ${\bf{\tilde u}}=(\tilde{u}_1,\ldots,\tilde{u}_6)$, ${\bf{\hat u}}=(\hat{u}_1,\ldots,\hat{u}_6)$
from the space
$\mathcal{C}(\overline{C_\varepsilon}\times[0,T]) \cap \mathcal{C}^{0,1}((C_\varepsilon \cup \Upsilon_{1,\varepsilon})\times[0,T]) \cap \mathcal{C}^{1,2}(C_\varepsilon\times(0,T))$ are called coupled
upper and lower solutions to the problem (\ref{start_prob}) if
${\bf \tilde u} \ge {\bf \hat u}$
(component-wise) and if they satisfy the differential inequalities
\begin{gather*}
\partial_t \tilde{u}_i - d_i \Delta \tilde{u}_i \ge f_i\big(t,x,\tilde{u}_i, [{\bf\tilde{u}}]_{a_i}, [{\bf\hat{u}}]_{b_i}\big),
\\
\partial_t \hat{u}_i - d_i \Delta \hat{u}_i \le f_i\big(t,x,\hat{u}_i, [{\bf\hat{u}}]_{a_i}, [{\bf\tilde{u}}]_{b_i}\big),
\end{gather*}
the boundary inequalities
\begin{gather*}
 d_i \partial_{\nu} \tilde{u}_i \ge g_i\big(x,\tilde{u}_i, [{\bf\tilde{u}}]_{c_i}, [{\bf\hat{u}}]_{d_i}\big),
\\
 d_i \partial_{\nu} \hat{u}_i \le g_i\big(x,\hat{u}_i, [{\bf\hat{u}}]_{c_i}, [{\bf\tilde{u}}]_{d_i}\big),
\end{gather*}
and the initial inequalities $\tilde{u}_i(0,x) \ge  0 \ge \hat{u}_i(0,x),$ and
$\tilde{u}_i = \hat{u}_i = 0$ on $\Gamma^\varepsilon_{0,\ell}\times(0,T)$ for every $i=1,\ldots,6.$
\end{definition}

Since  $f_i\big(t,x, 0, [{\bf 0}]_{a_i}, [{\bf\tilde{u}}]_{b_i}\big) \ge 0$ \ if $[{\bf\tilde{u}}]_{b_i}\ge [{\bf 0}]_{b_i}$ \ and
$g_i\big(x, 0, [{\bf 0}]_{c_i}, [{\bf\tilde{u}}]_{d_i}\big) \ge  0$ if $[{\bf\tilde{u}}]_{d_i}\ge [{\bf 0}]_{d_i}$ for $i=1,\ldots,6,$
\ ${\bf \hat u}\equiv \bf{0}$ is the lower solution to the problem (\ref{start_prob}).

Then the requirement for an upper solution is reduced to inequalities ${\bf \tilde u}(x,0) \ge {\bf 0},$
\begin{gather*}
\partial_t \tilde{u}_i - d_i \Delta \tilde{u}_i \ge f_i\big(t,x,\tilde{u}_i, [{\bf\tilde{u}}]_{a_i}, [{\bf 0}]_{b_i}\big),
\\
 d_i \partial_{\nu} \tilde{u}_i \ge g_i\big(x,\tilde{u}_i, [{\bf \tilde{u}}]_{c_i}, [{\bf 0} ]_{d_i}\big), \quad i=1,\ldots,6.
\end{gather*}
Obviously, that the solution to the linear system
\begin{equation}\label{upper_prob}
\left\{
    \begin{array}{rcll}
    \partial_t {\bf \tilde u} - \mathfrak{D} \Delta_x {\bf \tilde u}  &=& {\bf \tilde F}({\bf \tilde u})
                    &\mbox{in} \ C_{\varepsilon}\times (0,T),
                    \\[2mm]
                    \mathfrak{D}  \partial_{\nu} {\bf \tilde u} &=&   {\bf \tilde{G}}({\bf \tilde u})
                    &\quad \mbox{on} \ \Upsilon_{1}\times (0,T),
                    \\[2mm]
 \mathfrak{D}  \partial_{\nu} {\bf\tilde u} &=&   {\bf 0}
                    &\quad \mbox{on} \ \Upsilon_{1+\varepsilon}\times (0,T)
                    \\[2mm]
   {\bf \tilde u} &=& {\bf 0}& \text{on} \  \Gamma^\varepsilon_{0, \ell}\times (0,T),
 \\[2mm]
    \left. {\bf \tilde u}\right|_{t=0} &=&{\bf 0} &\mbox{in} \   C_{\varepsilon},
   \end{array}\right.
\end{equation}
is an upper solution to problem (\ref{start_prob}). In (\ref{upper_prob})
$$
{\bf \tilde F}({\bf u}) = 
\left(
\begin{array}{c}
 \psi(x, t) 
\\
\mu_1\, \Theta(x,t) \, u_1
\\
 0 
\\
\mu_6  +  \mu_7  u_4   
\\
\mu_8 + \mu_9 \, u_5 
\\
 \mu_{11} 
\end{array}
\right),
\quad
{\bf \tilde G}({\bf u}) =
\left(
\begin{array}{c}
\left(\gamma_1  + \gamma_2\right)\varphi_1(x)
\\
 0
\\
(\gamma_3 + \frac{\gamma_4 }{\lambda_{10}} )\varphi_2(x)
\\
 0
\\
 0
\\
 0
\end{array}
\right),
$$
and it is easily seen that for each $i\in \{1,\ldots,6\}$ the corresponding components of the vector-functions
${\bf \tilde F}=(\widetilde{f}_1,\ldots,\widetilde{f}_6),$ ${\bf \tilde G}= (\widetilde{g}_1,\ldots,\widetilde{g}_6)$ satisfy the inequalities
$$
\widetilde{f}_i({\bf u}) \ge f_i\big(t,x,u_i, [{\bf u}]_{a_i}, [{\bf 0}]_{b_i}\big)
, \quad
\widetilde{g}_i({\bf u}) \ge g_i\big(x,u_i, [{\bf u}]_{c_i}, [{\bf 0} ]_{d_i}\big) \quad \forall \ {\bf u} \in \Bbb R^6_+.
$$

In fact, the system (\ref{upper_prob}) is split into a system with respect $\widetilde{u}_1$ and $\widetilde{u}_2,$ and
four differential equations. Using classical results for linear parabolic boundary-value problems, we can state that
there exists a unique nonnegative solution ${\bf \tilde u}$ to the problem (\ref{upper_prob}) (this means the existence of an upper solution of (\ref{start_prob}))
and
\begin{equation}\label{bound-1}
  \sup_{(x,t) \in C_\varepsilon \times (0,T)} \left| {\bf \tilde u}(x,t) \right|   \le K(T).
\end{equation}

With the help of (\ref{bound-1}) it is easy to prove that for every $i\in\{1,\ldots,6\}$ there exit positive constants $K_i$ such that for the corresponding
components of the vector-functions  ${\bf F}$ and ${\bf G}$ the following inequalities hold:
\begin{gather}
  \left| f_i(t,x,{\bf u}) - f_i(t,x,{\bf v})  \right| \le K_i | {\bf u} - {\bf v}| \label{bound-2}
  \\
  \left| g_i(t,x,{\bf u}) - g_i(t,x,{\bf v})  \right| \le K_i | {\bf u} - {\bf v}| \label{bound-3}
\end{gather}
for all ${\bf 0}\le {\bf u} \le {\bf \tilde u},$ ${\bf 0}\le  {\bf v} \le {\bf \tilde{u}}$ and $(x, t) \in \overline{C}_\varepsilon \times [0,T].$

Hence all conditions of Theorem 7.2 from Ref.~\cite{Pao}[Chapter 9] are satisfied and we can state that the problem (\ref{start_prob}) has a unique solution ${\bf u}$ and
\begin{equation}\label{bound-4}
{\bf 0}\le {\bf u} \le {\bf \tilde{u}} \quad (\text{component-wise}) \quad \text{in} \ \overline{C}_\varepsilon \times [0,T].
\end{equation}

\begin{remark}
The existence and uniqueness of the solution to problem (\ref{start_prob})
can be also obtained from Ref.~\cite{LadColUra}[Charpter VII].
\end{remark}

\section{Asymptotic analysis of problem (\ref{start_prob})}\label{Sec4}

In the boundary condition for the component $u_1$ (see (\ref{reaction-term})) the parameter $\gamma_1$ is equal to
$$
\dfrac{R \, \alpha_0}{d_{L} \, k_0},
$$
where $d_L$ is the diffusion coefficient of LDLs in the intima
(it is equal to $29.89 \, cm^2 \, day^{-1}$ (see Refs.~\cite{CobSheMax,HaoFried-14})), \
 $k_0$ is the concentration of LDLs in the blood (estimated $k_0 \in \big(7 \times 10^{-4}\, , \, 1.9\times 10^{-3}\big)\ g \,cm^{-3})$
 at which the penetration rate is half of $V_{max},$ \
$\alpha_0$ is the penetration rate of LDLs from the blood into the intima through the damaged endothelium $\omega$
(estimated range is $(10^{-7}, \ 3\times 10^{-4})\, g \,cm^{-2} day^{-1}),$ \ $R$ is the radius of a blood vessel.
The size of blood vessels varies from $4\times 10^{-4}\, cm$ (radius of capillaries) to $1.2\,cm$ (radius of the aorta).
Atherosclerosis usually affects arteries of large and medium caliber. Therefore, we regard $R\in (10^{-2}\, , \, 1) \,cm.$
As a result,
\begin{equation}\label{gamma_1}
  \gamma_1  \in \left(10^{-8} \, , \, 1.5\times10^{-2} \right)
\end{equation}
is  a small parameter that can be
compared with $\varepsilon$ (nondimensional  thickness of the intima).
According to Ref.~\cite{arterial-1}, the ratio between the wall width and the radius of
blood vessels is relatively constant and it lays between $0.15$ and $0.20$. Then, $\varepsilon$
varies from $0.05$ to $0.07$.
Hence, we can introduce new values
$$
\gamma_1 = \varepsilon^{\varrho_1} \eta_1,\quad \gamma_2 = \varepsilon^{\varrho_2} \eta_2;
$$
similarly in the boundary conditions for $u_3$ and $u_4$:
$$
\gamma_3 = \varepsilon^{\varrho_3} \eta_3,\quad \gamma_4 = \varepsilon^{\varrho_4} \eta_4,
\quad \gamma_5 = \varepsilon^{\varrho_5} \eta_5
$$
where $\varrho_i \ge 1,$ and parameters $\eta_1,\ldots,\eta_5$ can be involved in the study
of active biochemical phenomena on the vessel wall like
 the endothelium dysfunction, adhesion of LDLs and monocytes, and their penetration or efflux.

For further asymptotic analysis we assume that  $\varrho_1=\ldots=\varrho_5 =: \varrho \ge 1.$
This assumption is made only for the sake of simplicity and it will be clear how to proceed in general case from further calculations.

Also small parameters are appeared through the coefficients
$$
d_5 = \frac{d_{M_1}}{d_L} = \frac{6.47 \times 10^{-5} \, cm^2 day^{-1}}{29.89 \, cm^2 day^{-1}} = 2.16\times 10^{-6} =: \varepsilon^{\tau_1},
$$
$$
d_6 = \frac{d_{F}}{d_L} = \frac{8.64 \times 10^{-7} \, cm^2 day^{-1}}{29.89 \, cm^2 day^{-1}} = 2.89 \times 10^{-8}=: \varepsilon^{\tau_2}.
$$
Because of $\varepsilon\in (0.05 , 0.07),$ we regard that $3 \le \tau_1 < \tau_2 < 2 \tau_1.$
Thus, we get the following parabolic system perturbed by the small parameter~$\varepsilon:$
\begin{equation}\label{start_prob-eps}
\left\{
    \begin{array}{rcll}
    \partial_t {\bf u}^\varepsilon - \mathfrak{D}_\varepsilon \, \Delta_x {\bf u}^\varepsilon  &=& {\bf F}({\bf u}^\varepsilon)
                    &\quad \mbox{in} \ C_{{\varepsilon}}\times (0,T),
                    \\[2mm]
 \mathfrak{D}_\varepsilon \, \partial_{\nu} {\bf u} &=&    {\varepsilon^\varrho} {\bf G}({\bf u}^\varepsilon)
                    &\quad \mbox{on} \ \Upsilon_{1}\times (0,T),
                    \\[2mm]
 \mathfrak{D}_\varepsilon \,  \partial_{\nu} {\bf u}^\varepsilon &=&   {\bf 0}
                    &\quad \mbox{on} \ \Upsilon_{1+\varepsilon}\times (0,T),
                    \\[2mm]
  {\bf u}^\varepsilon &=& {\bf 0}& \quad \text{on} \  \Gamma^\varepsilon_{0, \ell}\times (0,T),
 \\[2mm]
    \left. {\bf u}^\varepsilon\right|_{t=0} &=&{\bf 0} &\quad \mbox{in} \   C_{{\varepsilon}},
   \end{array}\right.
\end{equation}
where  $\mathfrak{D}_\varepsilon = {\rm diag}\left(1, d_2, d_3, d_4, {\varepsilon^{\tau_1}}, {\varepsilon^{\tau_2}}\right),$
${\bf G}({\bf u})= \big(g_1({\bf u}), 0, g_3({\bf u}), g_4({\bf u}), 0, 0\big)$  and
\begin{gather}\label{g-1}
g_1({\bf u}) = \varphi_1(x) \left(
 \eta_1 \dfrac{u_1}{1 + u_1} + \, \eta_2 \dfrac{u_6^{p_3}}{1 + u_6^{p_3}}\right),
\\\label{g-2}
g_3({\bf u}) = \varphi_2(x) \left(  \eta_3
\dfrac{u_2}{1 + u_2} + \, \eta_4 \dfrac{u_5^{p_1}}{(\lambda_{9} + u_4)(1 + u_5^{p_1})}\right),
\\ \label{g-3}
g_4({\bf u}) =
- \varphi_2(x)\,  \eta_5 \dfrac{  u_4^{p_2}}{1 + u_4^{p_2}}.
\end{gather}

{\sf Next our aim is to
construct the asymptotic approximation for the solution to the problem  (\ref{start_prob-eps}) as the small parameter $\varepsilon \to 0,$
\ to derive the limit problem $(\varepsilon =0),$ and to prove the corresponding asymptotic estimates.}

It should be noted that the limit process in the thin tube domain $C_\varepsilon$ is accompanied by the perturbed coefficients both in the differential equations and in the boundary conditions.

\subsection{Construction of the asymptotic approximation}

At first we rewrite the system (\ref{start_prob-eps}) in the cylindrical coordinates $(x_1, r, \vartheta):$
\begin{equation}\label{start_prob-eps-c}
\left\{
    \begin{array}{rcll}
    \partial_t {\bf u}^\varepsilon - \mathfrak{D}_\varepsilon \Delta_{x_1,r,\vartheta}{\bf u}^\varepsilon  &=& {\bf F}({\bf u}^\varepsilon)
                    &\mbox{in} \ \mathfrak{C}_{\varepsilon}\times (0,T),
                    \\[2mm]
 - \mathfrak{D}_\varepsilon  \partial_{r} {\bf u}^\varepsilon &=&  \varepsilon^\varrho {\bf G}({\bf u}^\varepsilon)
                    &\mbox{on} \ {\it \Upsilon}_{1}\times (0,T),
                    \\[2mm]
  \mathfrak{D}_\varepsilon  \partial_{r} {\bf u}^\varepsilon &=&  0
                    &\mbox{on} \ {\it \Upsilon}_{1+\varepsilon}\times (0,T),
 \\[2mm]
  {\bf u}^\varepsilon &=& {\bf 0}& \text{on} \  {\it \Gamma}^\varepsilon_{0, \ell}\times (0,T),
 \\[2mm]
   \partial_{\vartheta}^k{\bf u}^\varepsilon|_{\vartheta=0} &=& \partial_{\vartheta}^k{\bf u}^\varepsilon|_{\vartheta=2\pi}, & \ k =0, 1,
 \\[2mm]
    \left. {\bf u}^\varepsilon\right|_{t=0} &=&{\bf 0} &\mbox{in} \   \mathfrak{C}_{\varepsilon}.
   \end{array}\right.
\end{equation}
Here
$$
\Delta_{x_1,r,\vartheta} u =
\frac{\partial^2 u}{\partial x_1^2} + \frac{\partial^2 u}{\partial r^2} +  \frac{1}{r}\frac{\partial u}{\partial r}
    + \frac{1}{r^2} \frac{\partial^2 u}{\partial \vartheta ^2}
$$
is the Laplace operator in the cylindrical coordinates, the thin tube domain $C_\varepsilon$ is transformed into the thin
plate
$$
\mathfrak{C}_{\varepsilon} = \Big\{(x_1, r, \vartheta): \ \ x_1 \in (0, \ell ) , \quad r \in (1, 1+\varepsilon), \quad \vartheta \in (0, 2 \pi) \Big\},
$$
with the following parts of the boundary: ${\it \Gamma}^\varepsilon_{0,\ell}:= \partial \mathfrak{C}_{\varepsilon} \cap \big(\{ x_1=0\}\cup \{ x_1=\ell\}\big),$
\begin{gather*}
    {\it \Upsilon}_1= \big\{(x_1, \vartheta): \ \ x_1 \in (0, \ell ) , \quad r =1,  \quad \vartheta \in (0, 2 \pi) \big\}, \quad \text{and} \ \
   {\it \Upsilon}_{1+\varepsilon}.
\end{gather*}
Two cylindrical smooth surfaces $\omega_1$ and $\Omega_1$ are transfigured into smooth plane surfaces ${\it \omega}_1$ and ${\it \Omega}_1$ that belong to ${\it \Upsilon}_1.$

Using  the asymptotic approach to construct approximations for solutions of boundary-value problems in thin domains (see e.g. Refs.~\cite{Gol’den,Melnyk-2015,Mel-Klev}),
we propose the following approximation for the solution to problem (\ref{start_prob-eps-c}):
\begin{equation} \label{approx}
  {\bf u}^\varepsilon \approx {\bf R}^\varepsilon := {\bf v}(x_1, \vartheta, t) + {\varepsilon^2} \, {\bf W}(x_1,\tfrac{r-1}{{\varepsilon}},\vartheta,t),
\end{equation}
where ${\bf v}=(v_1,\ldots,v_6)$ and ${\bf W}=(W_1,\ldots,W_4,0,0).$

Let us substitute ${\bf R}^\varepsilon$  into (\ref{start_prob-eps-c}) instead of ${\bf u}^\varepsilon.$
Since the operator $\Delta_{x_1,r,\vartheta}$
in the variables $(x_1, \xi_2, \vartheta)$ takes the form
\begin{equation}\label{laplace}
  \Delta_{x_1,\xi_2,\vartheta} u = \frac{\partial^2 u}{\partial x_1^2} + \frac{1}{\varepsilon^2} \frac{\partial^2 u}{\partial \xi_2^2} +   \frac{1}{\varepsilon\xi_2+1} \frac{1}{\varepsilon} \frac{\partial u}{\partial \xi_2}
    + \frac{1}{(\varepsilon\xi_2+1)^2} \frac{\partial^2 u}{\partial \vartheta ^2},
\end{equation}
where $\xi_2= \frac{r-1}{\varepsilon},$ the collection of coefficients at the same power of $\varepsilon$ with regards (\ref{bound-2}) and (\ref{bound-3}) gives  the following four $(i=1,\ldots,4)$ one-dimensional Neumann boundary-value problems with respect to
$\xi_2:$
\begin{equation}\label{Neumann-prob-1}
\left\{
    \begin{array}{rcl}
- d_i \,\partial_{\xi_2\xi_2}^2 {W}_i(x_1,\xi_2,\vartheta,t) &=& f_i({\bf v}) - \partial_t v_i +  d_i \Delta_{x_1, \vartheta}{v}_i,   \ \ \ \xi_2\in (0, 1),
\\[1mm]
 - \,d_i \,\partial_{\xi_2} {W}_i|_{\xi_2=0} &=& {\delta_{1, \varrho}}\, g_i({\bf v}),\qquad \partial_{\xi_2} {W}_i|_{\xi_2=1} \, =\, 0,
\end{array}\right.
\end{equation}
and two ordinary differential equations
\begin{equation}\label{ord-diff-equ}
  \partial_t v_5 = f_5({\bf v}) \quad \text{and} \quad  \partial_t v_6 = f_6({\bf v}).
\end{equation}
In (\ref{Neumann-prob-1}), $\delta_{1, \varrho}$ is Kronecker's symbol (recall that $\varrho \ge 1$),
$
\Delta_{x_1, \vartheta} v = \frac{\partial^2 v}{\partial x_1^2} +  \frac{\partial^2 v}{\partial \vartheta ^2},
$
and the variables $x_1, \vartheta, t$ are regarded as parameters.

The solvability condition for the problem (\ref{Neumann-prob-1}) at a fixed index $i\in \{1,\ldots,4\}$ is given by the differential equation
 \begin{equation}\label{limit-eqs}
   \partial_t {v}_i - d_i \Delta_{x_1, \vartheta}v_i  = f_i({\bf v}) + \delta_{1, \varrho} \, g_i({\bf v}).
 \end{equation}

 Let $\{v_i\}_{i=1}^6$ be solutions of the respective differential equations (boundary and initial  conditions for them
will be determined later). Then  solutions of the problems (\ref{Neumann-prob-1}) exist and the additional relations
 \begin{equation}\label{midle-val}
\int_0^1 {W}_i(\xi_2,\cdot)\, d\xi_2 = 0, \quad  i=1,\ldots,4,
\end{equation}
 supply their  uniqueness. Since $g_2 \equiv 0,$ the solution ${W}_2 \equiv 0.$
Thus,  ${\bf W}=(W_1, 0, W_3, W_4, 0, 0).$

Equations  (\ref{limit-eqs}) and (\ref{ord-diff-equ}) form the limit system that is supplied by corresponding boundary and initial conditions. As a result, we get the coupled parabolic-ordinary system
 \begin{equation}\label{limit_prob}
\left\{
    \begin{array}{rcll}
    \partial_t {v}_i - d_i \Delta_{x_1,\vartheta}v_i  &=& f_i^{(0)}({\bf v}),
                    & \text{in}\  (0,\ell)\times(0, 2\pi)\times (0,T),
 \\[2mm]
   \partial_{\vartheta}^k v_i(x_1,\vartheta,t)|_{\vartheta=0} &=& \partial_{\vartheta}^k v_i(x_1,\vartheta,t)|_{\vartheta=2\pi},& k =0, 1,
 \\[2mm]
   v_i(0,\vartheta,t) &=& v_i(\ell,\vartheta,t)\ =\ 0, &  i=1,\ldots,4,
                    \\[2mm]
 \partial_t v_5 &=&   f_5({\bf v}), \quad  \partial_t v_6 = f_6({\bf v})
 \\[2mm]
    \left. {\bf v}\right|_{t=0} &=&{\bf 0} &\mbox{in} \  (0,\ell)\times(0, 2\pi),
   \end{array}\right.
\end{equation}
with the reaction terms
$$
{\bf F}^{(0)}= (f_1^{(0)},\ldots,f_4^{(0)}, f_5, f_6)
$$
$$
=
\left(
\begin{array}{l}
 - \, \Theta(x_1,1,\vartheta,t) \, v_1 +  \psi(x_1,1,\vartheta,t) + \varphi_1(x_1,\vartheta) \left(
  \delta_{1, \varrho_1} \eta_1 \dfrac{v_1}{1 + v_1} + \delta_{1, \varrho_2} \eta_2  \dfrac{v_6^{p_2}}{1 + v_6^{p_2}}\right)
\\[3mm]
 \mu_1\, \Theta(x_1,1,\vartheta,t) \, v_1 -  \mu_2 \dfrac{ u_2}{\lambda_1 + u_2} \, u_5 -  \mu_3 \dfrac{ u_2}{\lambda_2 + u_2} \, u_4
\\[3mm]
-  \mu_4 \, \dfrac{ v_3}{1 + v_3}   - \mu_5 \, \dfrac{ v_3}{\lambda_3 + v_3} -  \delta_3 \, v_3
 +  \varphi_2 \left(
\dfrac{\delta_{1, \varrho_3} \, \eta_3 \, v_2}{1 + v_2} +   \dfrac{\delta_{1, \varrho_4} \, \eta_4 \, v_5^{p_1}}{(\lambda_{9} + v_4)(1 + v_5^{p_1})}\right)
\\[3mm]
\mu_6 \dfrac{ v_3}{1 + v_3}  + \mu_7 \dfrac{ v_2}{\lambda_4 + v_2} \, v_4  -  \delta_4 \, v_4
 - \varphi_2(x_1,\vartheta) \delta_{1, \varrho_5}\, \eta_5 \, \dfrac{ v_4^{p_0}}{1 + v_4^{p_0}}
\\[3mm]
\mu_8 \dfrac{v_3}{\lambda_5 + v_3}  +  \mu_9 \dfrac{ v_2}{\lambda_6 + v_2} \, v_5 - \mu_{10} \dfrac{ v_5}{\lambda_7 + v_5}  - \delta_5 \, v_5
\\[3mm]
 \mu_{11} \dfrac{ v_5}{\lambda_8 + v_5}  - \delta_6 \, v_6
\end{array}
\right).
$$
Here we have written down the reaction terms of the limit problem without the assumption that $\varrho_1=\ldots=\varrho_5=:\varrho.$

The method of upper and lower solutions can be used to the system  (\ref{limit_prob}) as well (see Ref.~\cite{Pao}[\S~9.7]).
Similar as in subsection~\ref{subs-exist}, we can verify that the reaction terms ${\bf F}^{(0)}$ are quasi-monotone, ${\bf \hat{v}}\equiv {\bf 0}$ is a lower solution of (\ref{limit_prob}), and there exists an upper solution ${\bf \tilde{v}}.$ Then it follows from Theorem 7.3 (Ref.~\cite{Pao}[Chapter 9]) that the limit problem (\ref{limit_prob}) has a unique solution ${\bf v}$ and
\begin{equation}\label{bound-5}
{\bf 0}\le {\bf v} \le {\bf \tilde{v}} \quad (\text{component-wise}) \quad \text{in} \ [0,\ell]\times[0, 2\pi]\times [0,T].
\end{equation}
From positivity lemma (see~Ref.\cite{Pao}[Ch. 2]) it follows that either ${\bf v} >  {\bf 0}$ (component-wise) in $(0,\ell)\times [0, 2\pi]\times (0,T)$
or ${\bf v} \equiv  {\bf 0}$ that is impossible.

To be specific we will carry out further investigation for the more interesting case $\varrho =1$
(if $\varrho >1,$ then ${\bf W}\equiv 0).$

Thanks to the boundary conditions for the functions $v_1,\ldots,v_4,$ the initial value problem for $v_5$ at $x_1=0$ has the view
$$
\partial_t v_5(0,\vartheta,t)  =
- \mu_{10} \dfrac{ v_5(0,\vartheta,t)}{\lambda_7 + v_5(0,\vartheta,t)}  - \delta_5 \, v_5(0,\vartheta,t),\quad v_5(0,\vartheta,t)|_{t=0}=0.
$$
Due to Cauchy-Lipschitz theorem on existence and uniqueness of solutions to first-order equations with given initial conditions,
$v_5(0,\vartheta,t)=0$ for $\vartheta\in[0, 2\pi]$ and $t\in [0, T].$ Similarly, we verify that $v_6$ and $v_5$
vanish at $x_1=0$ and $x_1= \ell$ as well. Thus, ${\bf R}^\varepsilon|_{x_1=0} = {\bf R}^\varepsilon|_{x_1=\ell}.$

The approximation function ${\bf R}^\varepsilon$ takes into account the inhomogeneity
of the right-hand side ${\bf F}$ and the boundary conditions on ${\it \omega}_1$ and ${\it \Omega}_1$
that belong with their closures in the horizontal side ${\it \Upsilon}_1$ of the thin domain $\mathfrak{C}_\varepsilon.$ Since these boundary conditions are localized inside of ${\it \Upsilon}_1,$ the functions $W_1,  W_3, W_4$ vanish outside of $\overline{{\it \omega}_1}$ and
$\overline{{\it \Omega}_1}$ respectively. Therefore,
$\partial_r \big(v_1 + \varepsilon^2 W_1\big)|_{r=1} =0$ on ${\it \Upsilon}_1\setminus {\it \omega}_1,$ and
$\partial_r {\bf R}^\varepsilon|_{r=1} =0$ on ${\it \Upsilon}_1\setminus {\it \Omega}_1.$
Clearly that
$\partial_r {\bf R}^\varepsilon|_{r=1+\varepsilon} =0$
and
 $\partial_{\vartheta}^k{\bf R}^\varepsilon|_{\vartheta=0} =
\partial_{\vartheta}^k{\bf R}^\varepsilon|_{\vartheta=2\pi}$ for $k =0, 1.$

Now due to (\ref{limit-eqs}) we can rewrite problems for  $W_1,  W_3, W_4$ as follows
 \begin{equation}\label{Neumann-prob-2}
\left\{
    \begin{array}{rcl}
- d_i \,\partial_{\xi_2\xi_2}^2 {W}_i(x_1, \xi_2, \vartheta, t) &=& - g_i({\bf v}),   \ \ \  \xi_2\in (0, 1),
\\[1mm]
 - \,d_i \,\partial_{\xi_2} {W}_i|_{\xi_2=0} &=&  g_i({\bf v}), \qquad \partial_{\xi_2} {W}_i|_{\xi_2=1} \, =\,  0,
\\[1mm]
 \int_0^1 {W}_i\, d\xi_2&=&0.
\end{array}\right.
\end{equation}
Since $g_i({\bf v})|_{t=0}=0,$  we have $W_i|_{t=0}=0$ at $i\in\{1, 3, 4\}.$ Therefore,
${\bf R}^\varepsilon|_{t=0} = {\bf 0}.$

\subsection{Justification}

In virtue that $W_1,  W_3, W_4$ vanish in neighbourhoods of the sides $x_1=0$ and $x_1=\ell$
and ${\bf v} > {\bf 0}$ (see before), we can regard that
${\bf 0} \le {\bf R}^\varepsilon \le {\bf \tilde{u}}$ for $\varepsilon$ small enough.
Therefore, due to (\ref{bound-2}) and (\ref{bound-3}) we have
\begin{equation}\label{bound-4}
\left| {\bf F}(t,x,{\bf R}^\varepsilon) - {\bf F}(t,x,{\bf v}) \right| \le \varepsilon^2 c_1,
  \quad
\left| {\bf G}(t,x,{\bf R}^\varepsilon) - {\bf G}(t,x,{\bf v}) \right| \le \varepsilon^2 c_2
\end{equation}
in $\overline{\mathfrak{C}}_\varepsilon \times [0,T].$

\begin{remark}
Here and in what follows all constants $\{c_i\}$ in inequalities are independent of the parameter~$\varepsilon.$
\end{remark}

Substituting ${\bf R}^\varepsilon$ in the equations of the system (\ref{start_prob-eps-c}) and the respective boundary conditions on ${\it \omega}_1$
and ${\it \Omega}_1,$  and thanks to (\ref{Neumann-prob-1}) and (\ref{limit_prob}), we find
\begin{equation}\label{residual-1}
\partial_t {\bf R}^\varepsilon - \mathfrak{D}_\varepsilon \Delta_{x_1,r,\vartheta}{\bf R}^\varepsilon  -  {\bf F}({\bf R}^\varepsilon)=
 {\bf \Psi}^\varepsilon
\end{equation}
in $\mathfrak{C}_\varepsilon \times (0, T),$
\begin{gather}\label{bc-1}
  - d_i\partial_r \big(v_i + \varepsilon^2 W_i\big)|_{r=1} - \varepsilon g_i({\bf R}^\varepsilon) = \varepsilon \big( g_i({\bf v}) - g_i({\bf R}^\varepsilon)\big), \quad i\in \{1, 3, 4\},
  \\ \label{bc-2}
  - \partial_r v_i|_{r=1} = 0, \quad i\in \{2, 5, 6\},
\end{gather}
respectively on ${\it \omega}_1$ and ${\it \Omega}_1.$ In (\ref{residual-1})
${\bf \Psi}^\varepsilon:= \big({\bf F}({\bf v}) - {\bf F}({\bf R}^\varepsilon)\big) +  {\bf \Theta}^\varepsilon$
and  the components of ${\bf \Theta}^\varepsilon$ are as follow
\begin{gather*}
{\Theta}^\varepsilon_i = - d_i\left( \left(\frac{1}{r^2} - 1\right) \partial^2_{\vartheta \vartheta} {v}_i + \left(\frac{\varepsilon}{\varepsilon \xi_2 + 1} \partial_{\xi_2} {W}_i +
\frac{\varepsilon^2}{(\varepsilon \xi_2 + 1)^2} \partial^2_{\vartheta \vartheta} {W}_i
\right)\Big|_{\xi_2= \frac{r-1}{\varepsilon}}\right)
\\[2mm]
 +  \delta_{1,i}\,  \partial_{r}\psi(x_1,1+\theta,\vartheta,t)(r-1) + \varepsilon^2 \partial_t {W}_i
\end{gather*}
for $i=1,\ldots,4$ and
$$
{\Theta}^\varepsilon_i = -\varepsilon^{\tau_{i-4}} \left( \partial^2_{x_1 x_1} {v}_i + \frac{1}{r^2} \partial^2_{\vartheta \vartheta} {v}_i \right), \quad i=5, 6.
$$
Taking into account that $r\in (1, 1+ \varepsilon)$ and $\Psi^\varepsilon_i=\Theta^\varepsilon_i$ for $i\in \{2, 5, 6\},$ there exist positive constants $c_3$ and $\varepsilon_0$ such that for all $\varepsilon\in (0, \varepsilon_0)$
\begin{gather}\label{est-1}
\sup_{(x,t) \in \mathfrak{C}_\varepsilon \times (0, T)} \left| {\Psi}^\varepsilon_i(x,t) \right| \le c_3 \varepsilon \quad \text{for} \ \ i\in \{1,\ldots,4\},
\\
\sup_{(x,t) \in \mathfrak{C}_\varepsilon \times (0, T)} \left| {\Psi}^\varepsilon_i(x,t) \right| \le c_3 \, \varepsilon^{\tau_{i-4}} \quad \text{for} \ \ i\in \{5, 6\}. \label{est-2}
\end{gather}

Thus, the difference between the approximation function ${\bf R}^\varepsilon$ and
the solution ${\bf u}^\varepsilon$ to the problem  (\ref{start_prob-eps-c})  satisfies  the following relations:
\begin{equation*}
\left\{
    \begin{array}{rcll}
    \partial_t \big({\bf R}^\varepsilon - {\bf u}^\varepsilon\big) -
    \mathfrak{D}_\varepsilon \Delta_{x_1,r,\vartheta}\big({\bf R}^\varepsilon - {\bf u}^\varepsilon\big)
     &=& \big({\bf F}({\bf R}^\varepsilon) - {\bf F}({\bf u}^\varepsilon) \big) +
    {\bf \Psi}^\varepsilon
                        &\mbox{in} \ \mathfrak{C}_{\varepsilon}\times (0,T),
                    \\[2mm]
 - \mathfrak{D}_\varepsilon  \partial_{r} \big({\bf R}^\varepsilon - {\bf u}^\varepsilon\big) &=&
 \varepsilon\big({\bf G}({\bf R}^\varepsilon) - {\bf G}({\bf u}^\varepsilon)\big) +
  \varepsilon
 {\bf \Phi}^\varepsilon
                    &\mbox{on} \ {\it \Upsilon}_{1}\times (0,T),
                    \\[2mm]
  \mathfrak{D}_\varepsilon  \partial_{r} \big({\bf R}^\varepsilon - {\bf u}^\varepsilon\big) &=&  0
                    &\mbox{on} \ {\it \Upsilon}_{1+\varepsilon}\times (0,T),
 \\[2mm]
  \big({\bf R}^\varepsilon - {\bf u}^\varepsilon\big) &=& {\bf 0}& \text{on} \  {\it \Gamma}^\varepsilon_{0, \ell}\times (0,T),
 \\[2mm]
   \partial_{\vartheta}^k\big({\bf R}^\varepsilon - {\bf u}^\varepsilon\big)|_{\vartheta=0} &=& \partial_{\vartheta}^k\big({\bf R}^\varepsilon - {\bf u}^\varepsilon\big)|_{\vartheta=2\pi}, & \ k =0, 1,
 \\[2mm]
    \left. \big({\bf R}^\varepsilon - {\bf u}^\varepsilon\big)\right|_{t=0} &=&{\bf 0} &\mbox{in} \   \mathfrak{C}_{\varepsilon},
   \end{array}\right.
\end{equation*}
 or
\begin{equation}\label{difference}
\left\{
    \begin{array}{rcll}
    \partial_t\big({\bf R}^\varepsilon - {\bf u}^\varepsilon\big) - \mathfrak{D}_\varepsilon \, \Delta_x\big({\bf R}^\varepsilon - {\bf u}^\varepsilon\big)
      &=& \big({\bf F}({\bf R}^\varepsilon) - {\bf F}({\bf u}^\varepsilon) \big) +
    {\bf \Psi}^\varepsilon
                    &\mbox{in} \ C_{{\varepsilon}}\times (0,T),
                    \\[2mm]
 \mathfrak{D}_\varepsilon \, \partial_{\nu}\big({\bf R}^\varepsilon - {\bf u}^\varepsilon\big) &=&   \varepsilon\big({\bf G}({\bf R}^\varepsilon) -
 {\bf G}({\bf u}^\varepsilon)\big) + \varepsilon  {\bf \Phi}^\varepsilon
                    & \mbox{on} \ \Upsilon_{1}\times (0,T),
                    \\[2mm]
 \mathfrak{D}_\varepsilon \,  \partial_{\nu}\big({\bf R}^\varepsilon - {\bf u}^\varepsilon\big) &=&   {\bf 0}
                    & \mbox{on} \ \Upsilon_{1+\varepsilon}\times (0,T),
                    \\[2mm]
  {\bf R}^\varepsilon - {\bf u}^\varepsilon &=& {\bf 0}& \text{on} \  \Gamma^\varepsilon_{0, \ell}\times (0,T),
 \\[2mm]
    \left. \big({\bf R}^\varepsilon - {\bf u}^\varepsilon\big)\right|_{t=0} &=&{\bf 0} &\mbox{in} \   C_{{\varepsilon}},
   \end{array}\right.
\end{equation}
in the variables $x=(x_1,x_2,x_2).$ Here ${\bf \Phi}^\varepsilon:= {\bf G}({\bf v}) - {\bf G}({\bf R}^\varepsilon),$  ${\Phi}^\varepsilon_i=0, \ i=2, 5, 6,$ and
$$
\sup_{(x,t) \in \Upsilon_{1} \times (0, T)} \left| {\Phi}^\varepsilon_i(x,t) \right| \le c_4 \varepsilon^2, \quad i\in \{1, 3, 4\}.
$$

To deduce an estimate for the difference ${\bf R}^\varepsilon - {\bf u}^\varepsilon$ we use Theorem 7.3 in Chapter~5 of Ref.~\cite{LadColUra}, where an apriori estimate for the solution to the semilinear parabolic initial-boundary value problem
\begin{gather*}
\partial_t u - \sum\nolimits_{i,j=1}^na_{ij}(x,t,u)\partial^2_{x_i x_j}u = b(x,t,u,\nabla_x u) \quad \mbox{in} \ \Omega\times (0,T],
\\
\sum\nolimits_{i,j=1}^na_{ij}(x,t,u)\nu_j \partial_{x_i} u  = \psi_1(x,t,u) \quad \mbox{on} \ \partial\Omega\times (0,T],
\\
u|_{t=0}=\psi_0(x),
\end{gather*}
was proved under the following conditions for any $u$:
\begin{gather}\label{l-1}
0 \le \sum_{i,j=1}^na_{ij}(x,t,u) \xi_i \xi_j \le \mu_1 |{\bf \xi}|^2 \quad \text{for all}  \ (x,t)\in \overline{\Omega}\times (0,T],
\\\label{l-2}
u\, b(x,t,u, {\bf p}) \le \beta_0 |{\bf p}|^2 + \beta_1 |u|^2 + \beta_2 \quad \text{for all}  \ (x,t)\in \Omega\times (0,T],
\\\label{l-3}
\nu_1 |{\bf \xi}|^2 \le \sum_{i,j=1}^na_{ij}(x,t,u) \xi_i \xi_j, \quad
u \,\psi_1(x,t,u) \le \beta_3 |u|^2 + \beta_4 \quad \text{for all}  \ (x,t)\in \partial\Omega\times (0,T].
\end{gather}
Here $\nu_1, \mu_1 >0,$ $\beta_i\ge0, \ i\in\{0,1,\ldots,4\}.$ As a result, for any solution $u$ from $C^{2,1}(\Omega\times (0,T])\cap C(\overline{\Omega}\times [0,T])\cap C^{1,0}(\overline{\Omega}\times (0,T])$ the following estimate holds:
\begin{equation}\label{Lad-est}
  \max_{\overline{\Omega}\times [0,T]} |u(x,t)| \le \lambda_1 \exp(\lambda T) \max \left\{\sqrt{\beta_2}, \ \sqrt{\beta_4}, \ \max_{\overline{\Omega}} |\psi_0(x)| \right\},
\end{equation}
where the constants $\lambda_1$ and $\lambda$ are defined with the constants $\nu_1, \mu_1, \beta_0, \beta_1, \beta_3.$

With small changes the estimate carries over to every component of the solution ${\bf R}^\varepsilon - {\bf u}^\varepsilon$ to the problem (\ref{difference}).
For each $i\in \{1,\ldots,4\}$ the condition (\ref{l-1}) is satisfied. Then using  Young's inequality $(ab \le \varsigma a^2 + \tfrac{1}{4\varsigma}b^2$ for all $a, b \in \Bbb R$ and $\varsigma>0)$ with an suitable $\varsigma,$ (\ref{est-1}), (\ref{bound-4}),  properties of the components of (\ref{reaction-term-1}), (\ref{g-1})-(\ref{g-3})
and nonnegativity of ${\bf u}^\varepsilon$ and ${\bf R}^\varepsilon$,
we derive  analogues of the inequality (\ref{l-2})
\begin{gather}\label{in-10}
\big(R^\varepsilon_1 -  u^\varepsilon_1 \big) \big(f_1({\bf R}^\varepsilon) - f_1({\bf u}^\varepsilon) + \Psi^\varepsilon_1\big) \le \tfrac12 \big|R^\varepsilon_1 -  u^\varepsilon_1 \big|^2 + c_5 \varepsilon^2,
\\
\big(R^\varepsilon_2 -  u^\varepsilon_2 \big) \big(f_2({\bf R}^\varepsilon) - f_2({\bf u}^\varepsilon) + \Psi^\varepsilon_2\big) \le
C_2 \big|R^\varepsilon_2 -  u^\varepsilon_2 \big|^2 + c_5 \varepsilon^2 \notag
\\
 + \, \max_{\overline{\Omega}\times [0,T]}\big(\big|R^\varepsilon_1 -  u^\varepsilon_1 \big|^2+ \big|R^\varepsilon_4 -  u^\varepsilon_4 \big|^2 + \big|R^\varepsilon_5 -  u^\varepsilon_5 \big|^2 \big),
 \\
\big(R^\varepsilon_3 -  u^\varepsilon_3 \big) \big(f_3({\bf R}^\varepsilon) - f_3({\bf u}^\varepsilon) + \Psi^\varepsilon_3\big) \le
\tfrac12 \big|R^\varepsilon_3 -  u^\varepsilon_3 \big|^2 + c_5 \varepsilon^2,
\\
\big(R^\varepsilon_4 -  u^\varepsilon_4 \big) \big(f_4({\bf R}^\varepsilon) - f_4({\bf u}^\varepsilon) + \Psi^\varepsilon_4\big) \le
C(\varsigma_4)\big|R^\varepsilon_4 -  u^\varepsilon_4 \big|^2 + c_5 \varepsilon^2 \notag
\\
 + \, \varsigma_4\max_{\overline{\Omega}\times [0,T]}\big(\big|R^\varepsilon_2 -  u^\varepsilon_2 \big|^2 + \big|R^\varepsilon_3 -  u^\varepsilon_3 \big|^2 \big),
\end{gather}
and analogues of the second inequality in (\ref{l-3}) for $i\in\{1,3,4\}$
\begin{gather}
\big(R^\varepsilon_1 -  u^\varepsilon_1 \big) \varepsilon \big(g_1({\bf R}^\varepsilon) - g_1({\bf u}^\varepsilon) + \Phi^\varepsilon_1\big)\notag
\\
\le \varepsilon \Big(C_1 \big|R^\varepsilon_1 -  u^\varepsilon_1 \big|^2 + c_6 \varepsilon^4 + \max_{\overline{\Omega}\times [0,T]} \big|R^\varepsilon_6 -  u^\varepsilon_6 \big|^2\Big),
\\
\big(R^\varepsilon_3 -  u^\varepsilon_3 \big) \varepsilon \big(g_3({\bf R}^\varepsilon) - g_3({\bf u}^\varepsilon) + \Phi^\varepsilon_3\big) \le \varepsilon \Big(C_3 \big|R^\varepsilon_3 -  u^\varepsilon_3 \big|^2 + c_6 \varepsilon^4 \notag
\\
 + \, \max_{\overline{\Omega}\times [0,T]}\big(\big|R^\varepsilon_2 -  u^\varepsilon_2 \big|^2+ \big|R^\varepsilon_4 -  u^\varepsilon_4 \big|^2 + \big|R^\varepsilon_5 -  u^\varepsilon_5 \big|^2 \big)\Big),
\\
\big(R^\varepsilon_4 -  u^\varepsilon_4 \big) \varepsilon \big(g_4({\bf R}^\varepsilon) - g_4({\bf u}^\varepsilon) + \Phi^\varepsilon_4\big) \le
 \varepsilon \Big( \tfrac12 \big|R^\varepsilon_4 -  u^\varepsilon_4 \big|^2 + c_6 \varepsilon^4\Big).\label{in-11}
\end{gather}
Thus, due to (\ref{Lad-est}) we have
\begin{gather}\label{in-1}
  \max_{\overline{C}_\varepsilon\times [0,T]} \big|R^\varepsilon_1 -  u^\varepsilon_1 \big| \le
  C\Big( c_7 \varepsilon + \sqrt{\varepsilon} \max_{\overline{\Omega}\times [0,T]} \big|R^\varepsilon_6 -  u^\varepsilon_6 \big| \Big),
\\\label{in-2}
  \max_{\overline{C}_\varepsilon\times [0,T]} \big|R^\varepsilon_2 -  u^\varepsilon_2 \big| \le
C\Big( c_7 \varepsilon + \max_{\overline{\Omega}\times [0,T]}\big(\big|R^\varepsilon_1 -  u^\varepsilon_1 \big|+ \big|R^\varepsilon_4 -  u^\varepsilon_4 \big| + \big|R^\varepsilon_5 -  u^\varepsilon_5 \big| \big)\Big),
\\\label{in-3}
  \max_{\overline{C}_\varepsilon\times [0,T]} \big|R^\varepsilon_3 -  u^\varepsilon_3 \big| \le
C\Big( c_7 \varepsilon + \max_{\overline{\Omega}\times [0,T]}\big(\big|R^\varepsilon_2 -  u^\varepsilon_2 \big|+ \big|R^\varepsilon_4 -  u^\varepsilon_4 \big| + \big|R^\varepsilon_5 -  u^\varepsilon_5 \big| \big)\Big),
\\\label{in-4}
  \max_{\overline{C}_\varepsilon\times [0,T]} \big|R^\varepsilon_4 -  u^\varepsilon_4 \big| \le
  C\Big( c_7 \varepsilon
  + \sqrt{\varsigma_4} \max_{\overline{\Omega}\times [0,T]}\big(\big|R^\varepsilon_2 -  u^\varepsilon_2 \big|+ \big|R^\varepsilon_3 -  u^\varepsilon_3 \big| \big)\Big).
\end{gather}
For $i\in\{5, 6\}$ $v_i \equiv R^\varepsilon_i$ and $\Phi^\varepsilon_i\equiv 0,$ and we can repeat the proof of Theorem 7.3~\cite{LadColUra}
under conditions (\ref{l-1}) and (\ref{l-2}) only. In this case the constant $\lambda_1$ in (\ref{Lad-est}) is independent of $\varepsilon.$ As a result,
similarly as before we prove that
\begin{gather}
  \max_{\overline{C}_\varepsilon\times [0,T]} \big|v_5 -  u^\varepsilon_5 \big| \le
  C \Big( c_7 \varepsilon
 +\sqrt{\varsigma_5} \max_{\overline{\Omega}\times [0,T]}\big(\big|R^\varepsilon_2 -  u^\varepsilon_2 \big|+ \big|R^\varepsilon_3 -  u^\varepsilon_3 \big|\big)\Big),
  \\ \label{in-5}
\max_{\overline{C}_\varepsilon\times [0,T]} \big|v_6 -  u^\varepsilon_6 \big| \le
  C \Big( c_7 \varepsilon + \sqrt{\varsigma_6} \max_{\overline{\Omega}\times [0,T]}\big(\big|v_5 -  u^\varepsilon_5 \big|\Big).
 \end{gather}
Sequentially substituting inequality (\ref{in-1})--(\ref{in-5}) one in the other  and taking
$$
\varsigma_4= \frac{1}{4 \,C^4(C+2)^2},\quad \varsigma_5= \frac{1}{16 \,C^4(C+2)^2},\quad \varsigma_6= \frac{1}{4 \,C^4},
$$
we get
\begin{gather}\label{in-6}
   \max_{\overline{C}_\varepsilon\times [0,T]} \big|R^\varepsilon_4 -  u^\varepsilon_4 \big| \le
  C_4\Big( \varepsilon
  +  \max_{\overline{\Omega}\times [0,T]}\big(\big|v_5 -  u^\varepsilon_5 \big|+ \big|v_6 -  u^\varepsilon_6 \big| \big)\Big),
  \\ \label{in-7}
  \max_{\overline{C}_\varepsilon\times [0,T]} \big|v_5 -  u^\varepsilon_5 \big| \le
  C_5 \Big(  \varepsilon  + \max_{\overline{\Omega}\times [0,T]} \big|v_6 -  u^\varepsilon_6 \big|\Big)
  \\ \label{in-8}
\max_{\overline{C}_\varepsilon\times [0,T]} \big|v_6 -  u^\varepsilon_6 \big|  \le C_6 \, \varepsilon.
\end{gather}
Substituting (\ref{in-8}) in (\ref{in-1}), (\ref{in-7}) and then  obtained inequalities sequentially
in (\ref{in-6}), (\ref{in-2}) and (\ref{in-3}), we get  estimates for the other components.
Thus, we have proved the following theorem.
\begin{theorem}\label{mainTheorem}
There exist positive constants $C_0, \varepsilon_0$ such that for all values $\varepsilon\in (0, \varepsilon_0)$ the difference between the solution
${\bf u}^\varepsilon$ to the problem (\ref{start_prob-eps}) and the solution ${\bf v}$ to the
the limit coupled parabolic-ordinary system (\ref{limit_prob}) satisfies the following estimate
\begin{equation}\label{main-estmate}
 \max_{\overline{C}_\varepsilon\times [0,T]} \big|{\bf u}^\varepsilon - {\bf v} \big| \le C_0 \, \varepsilon.
\end{equation}
\end{theorem}

\begin{remark}
We think that the estimate (\ref{main-estmate}) could be obtained with the help of special integral representations using Green's functions for the corresponding linear initial-boundary value problems (see e.g. Ref.~\cite{Pao}[\S 9.6]). But the main difficulty will be
in establishing the dependence of the constant on the small parameter $\varepsilon$ in the proof of Lemma 9.6.1~\cite{Pao}.
\end{remark}

Multiplying each parabolic equation in (\ref{difference}) by the respective difference $R_i^\varepsilon - u^\varepsilon_i,$ integrating over $C_\varepsilon\times(0,T)$ by parts
and using the estimates (\ref{in-10})--(\ref{in-11}) and (\ref{main-estmate}), we find that for any $t\in (0, T)$
\begin{gather}
\frac12 \int_{C_{\varepsilon}} \big(R^\varepsilon_i(x,t) -  u^\varepsilon_i(x,t) \big)^2 dx
+
d_i \int_0^t \int_{C_{\varepsilon}} \big|\nabla(R^\varepsilon_i -  u^\varepsilon_i)\big|^2dx d\tau \notag
\\
  \le \int_0^T\int_{C_{\varepsilon}}\big| f_i({\bf R}^\varepsilon) - f_i({\bf u}^\varepsilon) + \Psi^\varepsilon_i \big|\,
  \big| R^\varepsilon_i -  u^\varepsilon_i \big| \, dx d\tau \notag
  \\+\,
  \varepsilon \int_0^T\int_{\Upsilon_1}\big| g_i({\bf R}^\varepsilon) - g_i({\bf u}^\varepsilon) + \Phi^\varepsilon_i \big|\, \big| R^\varepsilon_i -  u^\varepsilon_i\big|\,  d\sigma_x d\tau \le C \varepsilon^3, \quad i=1,\ldots4. \label{ine-1}
\end{gather}
Hence, the following estimate for gradients holds.
 \begin{corollary}\label{energy estimate}
 $$
\sum_{i=1}^4
\big\| \nabla u^\varepsilon_i -  \nabla R^\varepsilon_i  \big\|_{L^2(C_\varepsilon\times(0,T))}  \le C_1 \, \varepsilon^\frac32.
$$
\end{corollary}

\section{Conclusions}\label{Sec5}

Overview mathematical approaches  presented in Section~\ref{Sec2}  and the results obtained in this paper suggest that the development of new models of the atherosclerosis should be based on known scenarios of enzymatic reactions, namely on the corresponding differential equations describing these reactions. Furthermore, those mathematical models should be  nondimensionalised to introduce small (or large) scales in the system dynamics with a multiscale perspective to derive the respective simplified limit problem. In our case this is the small parameter $\varepsilon$ (see Section~\ref{Sec4}).

It is very important for a proposed multi-scale method to justify its stability and accuracy. The proof of the
error estimate between the constructed approximation and the exact solution is a general principle that has been
applied to the analysis of the efficiency of a multi-scale method. In the present paper, the asymptotic
approximation ${\bf R}^\varepsilon$ (see (\ref{approx})) for the solution ${\bf u}^\varepsilon$ to problem (\ref{start_prob-eps}) is constructed and justified. The results obtained in Theorem~\ref{mainTheorem} and Corollary \ref{energy estimate} argue that
 it is possible to replace the perturbed parabolic system (\ref{start_prob-eps}) in the thin tube domain $C_\varepsilon$ with
the corresponding limit problem (\ref{limit_prob})  in the rectangle $(0, \ell)\times (0,T)$ with the sufficient accuracy measured by the parameter $\varepsilon$ characterizing the non-dimensional  thickness of the intima layer.

From the viewpoint of mathematical modelling, this work is only an early start to the study of the
diverse roles of macrophages in atherosclerosis, because there remain a few important problems unexplored.
Among these open problems the following problems are most important and interesting  in author's opinion.
\begin{itemize}
  \item
From a practical point of view it is very important to stabilize the development of atherosclerosis.
A basic mathematical question about this problem is whether the solution ${\bf u}^\varepsilon$ to the problem (\ref{start_prob-eps}) converges
to the steady state as $t\to +\infty.$ Obviously, first of all this will depend on the parameters $\eta_1,\ldots,\eta_5$ and $p_1,$ $p_2,$ $p_3$
that are responsible for active biochemical phenomena on the vessel wall (the endothelium dysfunction, adhesion of LDLs and monocytes,
penetration of monocytes, and efflux of $M_2$-macrophages), and also on the other parameters in the problem (\ref{start_prob-eps}).

\smallskip

Due to the asymptotic results obtained in this paper,  the problem is reduced to the study of the time evolution of the solution to the limit problem (\ref{limit_prob}).  This will be a simpler problem, since there are no nonlinear boundary conditions. It is known (see Ref.~\cite{Pao}, Chapter 10) that in the case of a homogeneous Neumann boundary condition the asymptotic behavior (as $t\to +\infty)$ is compared with the behavior of the solution to the corresponding ordinary differential system.

\smallskip
\item
To study the influence of the parameters of the limit problem (\ref{limit_prob}) on the velocity of  the atherosclerosis development. In this connection, it will be interesting to examine the unboundedness of the solution to the problem (\ref{limit_prob}) with emphasis on its blowing-up behavior in finite time.

\smallskip
  \item
During the atherosclerosis development, the intimal thickness and the endothelial damage area grow and this should be taken into account.
We did this in the boundary condition (\ref{bc-4}) through the special  feedback term. However, it will be more appropriate to consider the corresponding free boundary problem modeling both the growth  of the  intimal tissue and damage area.
 For this the methods proposed in Refs.~\cite{Baz-Fried-1,Baz-Fried-2,Byrne-99,Bel-Prez,Bel_growth} will be helpful.
\end{itemize}

Of course, the introduced model is an incomplete representation of all the processes of the development of atherosclerosis. However, we believe that it reflects many of the main features of the  atherosclerosis and can give new understanding into the complex interactions of the disease.
In addition, we can complicate the model taking into account other features. For instance, when the effect of diffusion and convection are both taken into consideration then the new term ${\bf b}\cdot \nabla u_i$ appears in the corresponding differential equation (see Refs.~\cite{Banks,Fowler}).
In this case the method of upper and lower solutions can also apply to reaction diffusion-convection system with quasimonotone functions
(see Ref.~\cite{Pao}, Chapter 12). Also it would be interesting to include some chemotaxis terms in the model (see Refs.~\cite{Bel_chemotaxis,Jaeger_chemotaxis,Keller_Segel,Murray-2}).

\section*{Acknowledgment}
This research was begun at the University of Stuttgart under support of the Alexander von Humboldt Foundation in the summer of 2015.
The author  is grateful to Prof. Christian Rohde for the hospitality and wonderful working conditions. Then this research was continued
in the framework of the Marie Curie IRSES project "EU-Ukrainian Mathematicians for Life Sciences" (FP7-People-2011-IRSES Project
number 295164) and of the  Marie Curie RISE project "Approximation Methods for Modelling and Diagnosis Tools".
Some results were reported on the workshops of these projects.

\end{document}